  \documentclass{article}
 \usepackage{amssymb,amsmath,texdraw}

\begin{document}

 \def\Vert{{\rm Vert}}
 \def\Edge{{\rm Edge}}
 \def\N{{\cal N}}
 \def\s{{\rm symp}}
 \def\R{{\mathbb R}}
 \def\C{{\mathbb C}}
  \def\Z{{\mathbb Z}}
   \def\Q{{\mathbb Q}}
 \def\Spher{{\rm Hier}}
 \def\Abs{{\rm Abs}}
 \def\H{{\cal H}_\lambda}
 \def\E{{\cal E}_\lambda}
  \def\EE{{\frak E}_\lambda}
  \def\ov{\overline}
  \def\T{{\cal T}}
  \def\ge{\geqslant}
  \def\le{\leqslant}
  \def\cH{{\cal H}}
  \def\cK{{\cal K}}

 \newcounter{sec}
 \renewcommand{\theequation}{\arabic{sec}.\arabic{equation}}

 \def\1{   \medskip
  \begin{texdraw}

  \drawdim cm \linewd 0.03

 \move (0 2)   \fcir f:0 r:0.1
  \lpatt(0.1 0.1) \lvec (0 4)
     \fcir f:0 r:0.1  \rmove(0.3 0)    \htext{$a$} \rmove (-0.3 0)
  \lpatt( )  \lvec (0 6)
     \fcir f:0 r:0.1
    \move(0 4)
    \lvec (1 6)
     \fcir f:0 r:0.1
     \move(0 4)
     \lvec (-1 6)
      \fcir f:0 r:0.1
     \move(1 6) \rlvec ( 0.5 2)  \fcir f:0 r:0.1
     \move(1 6) \rlvec ( 0 2)  \fcir f:0 r:0.1
       \move(-1 6) \rlvec ( -0.5 2)  \fcir f:0 r:0.1
       \move(-1 6) \rlvec ( 0 2)  \fcir f:0 r:0.1
        \move(0 6) \rlvec ( 0 2)  \fcir f:0 r:0.1
          \move(0 6) \rlvec ( 0.5 2)  \fcir f:0 r:0.1
           \move(0 6) \rlvec ( -0.5 2)  \fcir f:0 r:0.1
          \move (-1 1) \htext{\sf A branch}

   \move (7 2)   \fcir f:0 r:0.1
    \lpatt(0.1 0.1) \lvec (7 4)
       \fcir f:0 r:0.1
       \move (7 4)\rlvec(1 -2)  \fcir f:0 r:0.1
       \move (7 4)\rlvec(-1 -2)  \fcir f:0 r:0.1
    \lpatt( )
    \move(7 4)        \rmove(0.3 0)    \htext{$a$} \rmove (-0.3 0)
    \lvec (7 6)
       \fcir f:0 r:0.1
      \move(7 4)
      \lvec (8 6)
       \fcir f:0 r:0.1
       \move(7 4)
       \lvec (6 6)
        \fcir f:0 r:0.1
       \move(8 6) \rlvec ( 0.5 2)  \fcir f:0 r:0.1
       \move(8 6) \rlvec ( 0 2)  \fcir f:0 r:0.1
         \move(6 6) \rlvec ( -0.5 2)  \fcir f:0 r:0.1
         \move(6 6) \rlvec ( 0 2)  \fcir f:0 r:0.1
          \move(7 6) \rlvec ( 0 2)  \fcir f:0 r:0.1
            \move(7 6) \rlvec ( 0.5 2)  \fcir f:0 r:0.1
             \move(7 6) \rlvec ( -0.5 2)  \fcir f:0 r:0.1

         \move(6.5 1) \htext{\sf A bush}

    \move(3.5 0) \htext{\sc Picture 1}

  \end{texdraw}\medskip}

\def\2{ \medskip
  \begin{texdraw}
    \drawdim cm \linewd 0.03
     \fcir f:0 r:0.1
    \move(0 0) \lvec(2 0)   \fcir f:0 r:0.1
    \lpatt (0.1 0.1)
    \move(0 0) \rlvec(-1 1)     \fcir f:0 r:0.1
    \move (2 0) \rlvec(1 1)      \fcir f:0 r:0.1
     \lpatt()
     \move(-1 1) \rlvec (0 2)    \fcir f:0 r:0.1
      \move(-1 1) \rlvec (-2 0)    \fcir f:0 r:0.1
        \move(3 1) \rlvec (0 2)    \fcir f:0 r:0.1
         \move(3 1) \rlvec (2 0)    \fcir f:0 r:0.1

         \move(0 0) \rlvec(-1 -1)     \fcir f:0 r:0.1
          \move(2 0) \rlvec(1 -1)     \fcir f:0 r:0.1
       \move(-1 -1) \rlvec (0 -2)    \fcir f:0 r:0.1
        \move(-1 -1) \rlvec (-2 0)    \fcir f:0 r:0.1
           \move(3 -1) \rlvec (0 -2)    \fcir f:0 r:0.1
            \move(3 -1) \rlvec (2 0)    \fcir f:0 r:0.1
      \linewd 0.005
      \move(-2 2) \lcir r:1.7
      \move(4 2)  \lcir r:1.7
       \lpatt (0.05 0.05)
       \move(1 2) \ravec (1.2 0)
        \move(1 2) \ravec (-1.2 0)

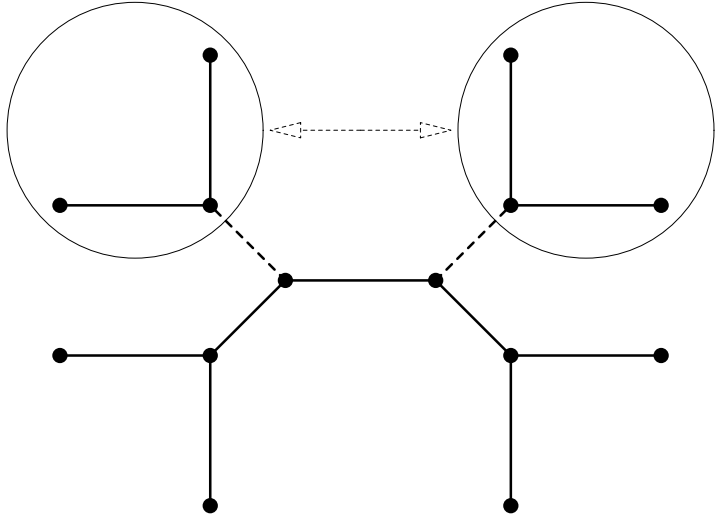
        \move(-4 -4) \htext{{\sc Picture 2.} \sf An example of
        hierarchomorphism:
        a re-glueing of two branches}
   \end{texdraw}\medskip}

   \def\3{ \medskip
    \begin{texdraw}
    \drawdim cm \linewd 0.01
    \lvec(4 0)
    \lvec (4 2.5)
    \move(4 0) \lvec(2.5 -1)
    \move (0 0) \lvec(2.5 -1)
    \move (0 0)   \lvec (4 2.5)
     \linewd 0.04
     \arrowheadtype t:F
     \move (0 0) \avec(2 0)
      \move (0 0) \avec(1.22 -0.5)
       \move (0 0) \avec (1.5 0.95)
     \move(1.22 1.2) \htext{$e_{b_1}$}
     \move(0.8 -0.9) \htext{$e_{b_2}$}
     \move(1.9 0.2) \htext{$e_a$}

     \move(2 -1.5) \htext{\sc Picture 3.}

   \end{texdraw}\medskip}

 \def\4{ \medskip
  \begin{texdraw}
   \drawdim cm \linewd 0.02
   \setunitscale  0.8
  \move(-7 0) \fcir f:0 r:0.1
  \rmove (0.3 0)\htext{$a$} \rmove (-0.3 0)
  \rlvec(0 2)  \fcir f:0 r:0.1
  \rmove (0.3 0)\htext{$b$} \rmove (-0.3 0)
  \rlvec(0 1.5)  \fcir f:0 r:0.1
  \rmove (0.3 0)\htext{$c$} \rmove (-0.3 0)
   \rlvec(0 1.5)  \fcir f:0 r:0.1
    \rmove (0.3 0)\htext{$d$} \rmove (-0.3 0)
   \rmove(0 -3) \rlvec(-1 0) \fcir f:0 r:0.1
   \rmove (-0.5 0)\htext{$h$} \rmove (0.3 0)
    %space pictures
    \linewd 0.005
    \move(0 0)
    \rlvec(-4 0)
    \rlvec(-3 -2)
     \rmove(-0.2 -0.5) \htext{$N_a$} \rmove(0.2 0.5)  \fcir f:0 r:0.1
    \rlvec(4 0)
       \rmove(-0.2 -0.5) \htext{$N_b$} \rmove(0.2 0.5)  \fcir f:0 r:0.1
    \rlvec(3 2) %one parallelogram
    \rlvec(0 3)
     \rmove(0.2 -0.2) \htext{$N_d$} \rmove(-0.2 +0.2)  \fcir f:0 r:0.1
    \rlvec(-4 0)
    \rlvec(0 -3)
    \rmove(4 0) %second parallelogram
    \rmove(0 3)
    \rlvec(-3 -2)
    \rmove(0.5 -0.2) \htext{$N_c$} \rmove(-0.5  0.2)  \fcir f:0 r:0.1
    \rlvec(0 -3)
    %second space picture

      \linewd 0.005
      \move(6.5 0)
        \rmove(0.2 -0.2) \htext{$N_b$} \rmove(-0.2 0.2)  \fcir f:0 r:0.1
      \rlvec(-4 0)
         \rmove(0.1 -0.5) \htext{$N_h$} \rmove(-0.1 0.5)  \fcir f:0 r:0.1
      \rlvec(-3 -2)
      \rlvec(4 0)
      \rmove(0.2 -0.2) \htext{$N_a$} \rmove(-0.2 0.2)  \fcir f:0 r:0.1
      \rlvec(3 2)

      %one parallelogram
      \rlvec(0 3)
        \rmove(0.2 -0.0) \htext{$N_c$} \rmove(-0.2 0.)  \fcir f:0 r:0.1
      \rlvec(-4 0)
      \rlvec(0 -3)
      \rmove(4 0) %second parallelogram
      \rmove(0 3)
      \rlvec(-3 -2)
      \rlvec(0 -3)
      %second space picture

      \move(-8 -3.5)\htext{\sc Picture 4.}
     \move(-8 -4)\htext{\sf Five points $N_a$, $N_b$, $N_c$, $N_d$,
     $N_h$
      span a 4-dimensional subspace in the affine Hilbert space $K$.}
       \move(-8 -4.5)\htext{\sf
      We portray the relative positions of $N_a$,  $N_b$, $N_c$, $N_d$
      in the corresponding 3 dimensional  space,}
       \move(-8 -5)\htext{\sf
      and also
      the relative positions  $N_a$, $N_b$, $N_c$,   $N_h$
      in (another) 3-dimensional space}

   \end{texdraw}\medskip}

   \def\6{ \medskip
     \begin{texdraw}
      \drawdim cm \linewd 0.03
      \fcir f:0 r:0.1
     {\lpatt (0.1 0.1)\lvec(4 0)}
       \fcir f:0 r:0.1
     \move(4 0)\lpatt()\lvec(6 0) \fcir f:0 r:0.1
      \move(4 0)\lpatt()\lvec(6 1) \fcir f:0 r:0.1
      \move(4 0)\lpatt()\lvec(6 -1) \fcir f:0 r:0.1
     \move(0 -0.5) \textref h:C v:T \htext{$\xi$}
      \move(4 -0.5) \textref h:C v:T \htext{$u$}
       \move(6 -1.5) \textref h:C v:T \htext{$v_k$}
       \lpatt (0.1 0.1)
       \move(6 0)\rlvec(2 0.5)
        \move(6 0)\rlvec(2 -0.5)
        \move(6 0)\rlvec(2  0)
        \move(6 1)\rlvec(2 0.5)
          \move(6 1)\rlvec(2 0)
         \move(6 -1)\rlvec(2 -0.5)
           \move(6 -1)\rlvec(2 0)

           \move(1 -2)\htext{\sc Picture 6.}

     \end{texdraw}\medskip }

     \def\5{\medskip
    \begin{texdraw}
    \drawdim cm
     \linewd 0.03
     \move(1 0) \rlvec(0 1) \rlvec(-1 0)
     \move(1 1) \rlvec(1 1) \rlvec(-1 1) \rlvec (-1 0)
     \move(1 3) \rlvec(0 1)
     \move(2 2) \rlvec(1 0)\rlvec(1 -1)
    \move(3 2) \rlvec(1 1)\rlvec(0 1)
     \move(4 3)\rlvec(1 0)\rlvec(1 -1)
     \move(5 3)\rlvec(1 1)

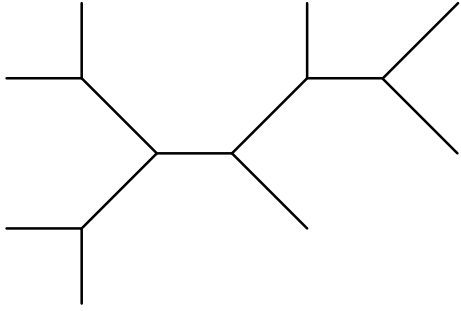
     \move(1 -1)\htext{\sc Picture 5. \sf A complete subtree in the dyadic
  Bruhat-Tits   tree $\T_2$}

    \end{texdraw}
    \medskip}

 \begin{center}

 {\Large Groups of hierarchomorphisms of trees
  and related Hilbert spaces}

 \medskip

{\large Yurii A. Neretin}\footnote%
{partially supported by the
grant
NWO 047-008-009}

 \end{center}

{\bf 0.1. Hierarchomorphisms (spheromorphisms).}
 The Bruhat--Tits tree  $\T_p$ is a infinite tree such that any
 vertex belongs to $(p+1)$ edges.
 As was observed by Cartier \cite{Car}, the groups
 $Aut(\T_p)$ of automorphisms of the trees $\T_p$
 are analogues of real and $p$-adic groups of rank 1
 (as ${\rm SL}_2(\R)$, ${\rm SL}_2(\C)$,
 ${\rm O}(1,n)$, ${\rm SL}_2(\Q_p)$ etc.).
  The representation theory of   $Aut(\T_p)$
  was developed in Cartier's \cite{Car} and Olshansky's
  \cite{Ols1} papers.
  In fact,
 the group $Aut(\T_p)$ is essentially simpler
 than the rank 1 groups on locally compact fields, but
 many nontrivial phenomena related to rank 1 groups survive
 for the group of automorphisms of Bruhat--Tits trees.

 The absolute of the Bruhat--Tits tree is an analogue of the boundaries
 of rank 1 symmetric spaces, in particular, the absolute
 is an analogue of the circle.
 The group of hierarchomorphisms\footnote{In \cite{Ner2},
 there was proposed the term {\it 'ball-morphisms'},
 which is difficult for pronouncement.
 In English translation, it was replaced by {\it 'spheromorphism'}.
 I want to propose the neologism
  {\it 'hierarchomorphism'}, this a map
  regarding hierarchy of balls on the absolute;
   see below
  Subsection 5.1.}
  $\Spher(\T_p)$
 (defined in \cite{Ner1}) is a tree analogue  of
 the group ${\rm Diff} (S^1)$ of diffeomorphisms of the circle.
 The group  $\Spher(\T_p)$
  consists of homeomorphisms of the absolute of $\T_p$
 that can be extended to the whole Bruhat--Tits tree except
 a finite subtree. It turns out to be (\cite{Ner1}, \cite{Ner2}), that
 the representation theory of ${\rm Diff} (S^1)$
 partially survives for the groups  $\Spher(\T_p)$.

 In fact, the group  $\Spher(\T_p)$  contains the
 group
 of locally analytic diffeomorphisms of $p$-adic line
 (see \cite{Ner2}),
 and this partially explains the similarity of ${\rm Diff} (S^1)$
  and $\Spher(\T_p)$.\footnote{%
  Another heuristic explanation can be obtain
  by the monstrous degeneration construction from \cite{Kapo},
  chapter 9;
  the Lobachevsky plane can be degenerated to
  the universal $\R$-tree.}

 The following  facts $1^\circ$-$4^\circ$
  are known about the groups $\Spher(\T_p)$.
 The phenomena $1^\circ$-$3^\circ$ are an exact  reflection
 of the representation theory of ${\rm Diff} (S^1)$,
 the last phenomenon now does not have a visible real analogue.

 \smallskip

$1^\circ$. (\cite{Ner1}, \cite{Ner2}) Denote by ${\rm O}(\infty)$
the group of all orthogonal operators
in a real Hilbert space $H$. Denote by ${\rm GLO}(\infty)$
the group of all invertible operators in $H$ having the form
$A=B+T$, where $B\in {\rm O}(\infty)$ and $T$ has finite rank.
Denote by $H_\C$ the complexification of $H$. Denote by
${\rm UO}(\infty)$ the group of all unitary operators in $H_\C$
having the form $A=B+T$, where $B\in {\rm O}(\infty)$ and $T$ has finite rank.
There exist some series of embeddings
$$   \Spher(\T_p)\to {\rm GLO}(\infty);\qquad
       \Spher(\T_p)\to {\rm UO}(\infty).$$
This allows to apply  the second quantization machinery
(see \cite{Ols3}, \cite{Ners}, \cite{Nerb})
for obtaining unitary representations
of   $\Spher(\T_p)$.

\smallskip

$2^\circ$. Embeddings
$\Spher(\T_p)\to {\rm GLO}(\infty)$ allow to develop
a theory of fractional diffusions with a Cantor set
time (the Cantor set appears as
the absolute of the tree).
I never wrote a text on this topic,
 but, on the whole, the picture
here is quite parallel to fractional diffusions
with real time (see, \cite{Ner3}).

\smallskip

$3^\circ$.(Kapoudjian, \cite{Kap1}, \cite{Kap4})
 There exists a $\Z/2\Z$-central extension of $\Spher(\T_p)$.

\smallskip

$4^\circ$.(Kapoudjian, \cite{Kap3})
 Consider the dyadic Bruhat--Tits tree $\T_2$.
There exists a canonical action of the group $\Spher(\T_2)$
on the inductive limit of the Deligne--Mumford \cite{DM} moduli
spaces $\lim_{n\to\infty}{\cal M}_{0,2^n}$ of $2^n$ point
configurations on the Riemann sphere. This construction
also has two  versions over $\Bbb R$. The first
variant is
an action on the
inductive limit of Stasheff associahedrons
(\cite{Sta}). The second variant is an action on the inductive limit
of the spaces constructed by Davis, Janiszkiewicz, Scott
(\cite{Dav}).
The last case is  most
interesting, since this real space has an interesting topology.

\smallskip

{\bf 0.2. The purposes of this paper.}
This paper has two purposes. The first aim  is to construct a new
series of embeddings of the groups of hierarchomorphisms
to the group  ${\rm GLO}(\infty)$. By
the Feldman-Hajek theorem
(see \cite{Sko}), this gives constructions of
unitary representations of groups of
hierarchomorphisms, but we do not discuss
this subject.

There exists the wide and nice theory of actions
of groups  on trees (see \cite{Ser}, \cite{Sha1}, \cite{Sha2},
\cite{Kapo}). It is clear that a hierarchomorphism
type extension can be constructed for any group $\Gamma$
acting on a tree (and even on an $\R$-tree), it is sufficient
to allow to cut a finite collection of edges.
The second purpose of this paper\footnote{see also
the recent preprint of Nekrashevich \cite{Nek}.}
 is to understand,
is  this "hierarchomorphization" of arbitrary group $\Gamma$
a reasonable object?

One example of such "hierarchomorphization"
is quite known, this is the Richard Thompson group \cite{Tho},
which firstly appeared as an counterexample
in theory of discrete groups.
Later it became clear, that this group is not
 a semipathological counterexample,
but a rich and unusual object (see works of Greenberg, Ghys,
Sergiescu, Penner, Freyd, Heller and others
 \cite{GhS}, \cite{GrS}, \cite{Pen}, \cite{FH},
 \cite{BG}
 see also \cite{Can}), relation of
 the hierarchomorphisms and the Thompson group
 was observed by Sergiescu).

If the group $\Gamma$ is discrete, then the corresponding
group of hierarchomorphisms is a discrete Thompson-like group.
If the group $\Gamma$ is locally compact, then the group
of hierarchomorphisms (see some examples in \cite{Ner2})
is an "infinite dimensional group" (or, better, "large group")
 similar to the group
of diffeomorphisms of the circle or diffeomorfisms
of $p$-adic line.

 \smallskip

 {\bf 0.3. The structure of the paper.}
 Sections 1-2 contain preliminary definitions and examples.

 In Section 3, we define the groups of hierarchomorphisms
 of tree (this definition can be adapted also for $\R$-tree,
 but nontrivial constructions of Sections 4--6 do not
 survive in this case).

 In Section 4, we discuss a family of
  Hilbert spaces $\H(J)$, where $0<\lambda<1$,
 associated with a tree $J$. The space $\H(J)$ contains
 the (nonorthogonal) basis $e_a$
  enumerated by vertices $a$ of the tree, and the inner products
  of the vectors $e_a$, $e_b$ are given by
  $$\langle e_a, e_b\rangle=
  \lambda^{\{\text{distance between $a$ and $b$}\}}$$
  We show that the group of hierarchomorphisms of $J$
  acts in $\H(J)$ by operators of the class ${\rm GLO}(\infty)$.

  In Section 5, for sufficiently large $\lambda$ we
  construct an operator of the 'restriction to
  the absolute' in the space $\H$.
  %This Section is an imitation of the work \cite{NO}.

  In Section 6, we discuss the action of the group of hierarchomorphisms
  in spaces of functions (distributions)
   on the absolute.

    The results of
  Sections 4--5 are 'new' for the groups of hierarchomorphisms
  of the Bruhat--Tits trees. The construction of Section 6
  in for Bruhat--Tits trees coincides with \cite{Ner1}.

  \smallskip

  {\bf Acknowledgment.} I am grateful to
  V.Sergiescu and C.Kapoudjian for meaningfull discussions.
  I thank the administration of the Erwin Schr\"odinger Institute
  (Wien)
  and Institute Fourier (Grenoble),
   where this work was done, for hospitality.

 \medskip

 {\large\bf 1. Notation and terminology}

  \addtocounter{sec}{1}
 \setcounter{equation}{0}

 \smallskip

 {\bf 1.1. Symplicial trees.}
  {\it A simplicial tree} $J$ is a connected graph without circuits.

   By $\Vert(J)$ we denote the set of vertices of $J$.

   By $\Edge(J)$ we denote the set of edges of $J$.

   We say that two vertices $a,b\in\Vert(J)$ are {\it adjacent},
   if they are connected by an edge. We denote this edge by
   $[a,b]$.

   We assume that the sets $\Vert(J)$, $\Edge(J)$
    are countable or finite.  A simplicial
    tree is {\it locally finite}
    if any vertex $a$ belongs to finitely many edges.
    We  admit non locally finite trees.

   A {\it way} in $J$ is a sequence of {\it distinct} vertices
   $$\dots, a_1, a_2, a_3, \dots$$
   such that $a_j$, $a_{j+1}$ are adjacent.
   A way can be finite, or infinite to one side, or infinite to the both
   sides.

   For vertices $a,b$ there exists a unique way $a_0=a, a_1,\dots, a_k=b$
   connecting $a$ and $b$. We say that $k$ is the
   {\it symplicial distance} between $a$ and $b$. We denote the simplicial
   metrics by
   $$ d_{\s}(a,b).$$

   A {\it subtree} $I\subset J$ is a connected subgraph in
   the tree $J$.

   The {\it boundary} $\partial I$
    of a subtree $I\subset J$ is the set of
   all $a\in \Vert(I)$ such that there exists an edge $[a,b]$
   with $b\notin\Vert(I)$.

   A subtree $I\subset J$ is {\it right}, if the number of
   edges $[a,b]\in\Edge(J)$ such that $a\in I$, $b\notin I$ is finite.

   A subtree $I\subset J$ is a {\it branch} if there is a unique
   edge $[a,b]\in\Edge(J)$ such that $a\in \Vert(I)$, $b\not\in\Vert(I)$,
   see Picture 1.
   The vertex $a$ is called a {\it root} of the branch.
   If we delete an edge of the tree $J$, then we obtain two branches.

   A subtree $I\subset J$ is a {\it bush} if its boundary contain
   only one point $a$ (a {\it root}) and number of edges
   $[a,b]\in\Edge(J)$ such that
   $b\not\in I$ is finite, see Picture 1.

   \smallskip

   {\sc Lemma 1.1.} a) {\it The intersection of a finite family
   of right subtrees is a right subtree.}

   \smallskip

   b) {\it For a right subtree $I\subset J$,
   there exists
    a finite  collection  of edges $\ell_1$, \dots, $\ell_k\in\Edge(I)$
    such that $I$ without  $\ell_1$, \dots, $\ell_k$ is a union
    of bushes.}

    \smallskip

   {\sc Proof.} The statement a) is obvious.

       The statement b). Let $a_1$, \dots, $a_k$ be the boundary
       of $I$. Let $L\subset I$
        be the minimal subtree containing the vertices  $a_1$, \dots, $a_k$.
        It is sufficient to delete all edges of $L$.
        \hfill $\square$

       \smallskip

       \1

  We say that a tree $J$ is {\it perfect}  if any vertex  of $J$
  belongs to $\ge 3$ edges.
  Obviously, perfect trees  are infinite.

         \smallskip

  {\bf 1.2. Actions of groups on simplicial trees.}
  A bijection $\Vert(J)\to\Vert(J)$ is an {\it automorphism}
  of a simplicial tree
  $J$ if   the images of adjacent vertices are adjacent vertices.

  An action of a group $\Gamma$ on a simplicial tree is an embedding
  of $\Gamma$ to its group of automorphisms.

    \smallskip

  {\bf 1.3. Absolute.}  The absolute $\Abs(J)$
  of a tree is the set of points of
  the tree on infinity. Let us give the formal definition.

  We say that a {\it ray} is an infinite way $a_1,a_2,\dots$.
  We say that rays $a_1,a_2,\dots$ and  $b_1,b_2,\dots$
  are equivalent if there exist $k$ and a sufficiently
  large $N$ such that $b_j=a_{j+k}$ far all $j\ge N$.

  A point of an absolute is a class of equivalent ways.

  \smallskip

  {\bf 1.4. Metric trees.} Let $J$ be a simplicial tree. Let us assign
  a positive number $\rho(a,b)$ to each edge $[a,b]$.
  Let $a$, $c$ be arbitrary vertices of $J$, let $a_0=a,a_1,\dots,a_k=c$
  be the way connecting $a$ and $c$.
 We assume
 $$\rho(a,c)=\sum_{j=1}^k \rho(a_{j-1},a_{j})$$
 Obviously, $\rho$ is a metric on $\Vert(J)$.
 We call by {\it metric trees} the
 {\it countable}  spaces $\Vert(J)$
 equipped with the metrics $\rho$.

Obviously, the edges  of $J$ can be reconstructed using
the metric $\rho$.
Hence we prefer to think
that the edges are present in a metric tree as an additional
(combinatorial) structure.

{\sc Remark.} We also can assume that lengths of all
edges is 1, and thus a simplicial tree is a partial
case of  metric trees.

 \smallskip

 {\sc Remark.} In literature, sometimes the term {\it metric tree}
 is used
 in the quite different sense (for $\R$-trees).

 \smallskip

 A metric tree $J$ is {\it locally finite} if
 it is locally finite as a simplicial tree and
 for any $a\in\Vert(J)$ and each $C>0$ the set
 of vertices $b$ satisfying $\rho(a,b)<C$
 is finite.

 \smallskip

 {\bf 1.5. Actions of groups on metric trees.}
 Let $J$ be a metric tree. A bijection $\Vert(J)\to\Vert(J)$
 is an {\it automorphism} of $J$ if it preserves the distance
 (hence it automatically preserves the sructure of
  simplicial tree).

 An action of group $\Gamma$ on a metric tree $J$ is an embedding
 of $\Gamma$ to the group of automorphisms of $J$.

 \medskip

 {\large\bf 2. Examples of actions of groups on trees.}

  \addtocounter{sec}{1}
 \setcounter{equation}{0}

 \smallskip

 The purpose of this Section is to give a collection
 of examples for abstract constructions given  in Sections 3-6
 (all these examples are standard).
 For algebraic and combinatorial
 theory of actions of groups on trees, see \cite{Ser},
 \cite{Sha1}, \cite{Sha2}.

 \smallskip

 {\bf 2.1. Bruhat--Tits trees.}
 The Bruhat--Tits tree ${\cal T}_p$
 is the tree, in which each vertex belongs to $(p+1)$
 edges. The group $Aut ({\cal T}_p)$  of automorphisms of $\T_p$
 is a locally compact group. This group is similar
 to rank 1 groups over $\R$ and over $p$-adic fields.
  The representation theory
 of  $Aut ({\cal T}_p)$ and related harmonic analysis
 are well understood, see \cite{Car}, \cite{Ols1},
 \cite{Fig1}, \cite{Fig2}.

 \smallskip

 {\bf 2.2. The tree $\T_\infty$.}
 We denote by $\T_\infty$ the simplicial tree,
 in which each  vertex belongs to a countable set
 of edges. At first sight, the group $Aut(\T_\infty)$
 seems pathological. Nevertheless, it is a useful
 object as one of the simplest examples of infinite-dimensional
 groups, see \cite{Ols2}, \cite{Nerb}.
 This group is an imitation of the group ${\rm O}(1,\infty)$.

 \smallskip

 {\bf 2.3. The tree of free group.}
 Denote by $F_2$ the free group with two generators $\alpha$, $\beta$.
 Vertices of the tree $J(F_2)$ are numerated  by  elements of
 the group $F_2$.
 Vertices $v_p$, $v_q$ are connected by an edge
 if
 $$p=q\alpha^{\pm1} \qquad\text{or}\qquad p=q\beta^{\pm1}.$$
 Obviously, $J(F_2)$ is the Bruhat--Tits  tree $\T_3$.
 The group $F_2$ acts on the tree $J(F_2)$
 by the transformations
 $$r:\quad v_p\mapsto v_{rp}, \qquad\text{where}\quad r\in F_2.$$

Fix $l_1$, $l_2>0$. Assign the length
 $l_1$ to any edge
 $[v_p,v_{p\alpha}]$,
 and the length
 $l_1$ to any edge
 $[v_p,v_{p\beta}]$.
 Thus we obtain a metric tree with an action
 of $F_2$.

 \smallskip

 {\bf 2.4. Another tree of free group.}   Let us contract
 all the edges of the type $[v_p,v_{p\alpha}]$ of the tree
 $J(F_2)$ described in 2.3.  Thus, we obtain the action of $F_2$
 on $\T_\infty$.

 \smallskip

 {\bf 2.5 Dyadic intervals.}
 Vertices $V_{u;n}$ of the tree $J_2(\R)$ are enumerated by
 segments in $\R$ having the form
 $$ S_{u;n}=\left[ \frac u{2^n},  \frac {u+1}{2^n}\right], \qquad
 \text{where $u\in \Z, n\in \Z$}.
 $$
 We connect $V_{u;n}$ and $V_{w;n-1}$ by an edge
 if  $S_{w;n-1}\supset V_{u;n}$.

 Obviously, we obtain the simplicial tree $\T_2$.

 \smallskip

 {\bf 2.6. Balls on $p$-adic line}.
 Denote by $\Q_p$ the field of $p$-adic numbers, denote by ${\Bbb Z}_p$
the $p$-adic integers.
 Denote by $B_{a, k}$ the ball
 $$ |z-a|\le p^{-k}.$$

 \smallskip

 {\sc Remark.} The radius $p^{-k}$ is determined by the ball.
 But $B_{a,k}=B_{c,k}$ for any $c\in B_{a,k}$.

 \smallskip

{\bf 2.7. Tree of lattices.}
 Consider the $p$-adic plane $\Q_p^2$
  equipped with a skew symmetric bilinear
 form $A(v,w)$. Denote by ${\rm Sp_2(\Q_p)}$ the group
 of linear transformations preserving
 the form $A(v,w)$.

 A {\it lattice} in  $\Q_p^2$ is a compact subset $R \subset\Q_p^2$
 having the form
 $$\Q_pv \oplus \Q_pw; \qquad\text{where $v,w$ are not proportional}.$$
 We say that a lattice $R$ is {\it self-dual}
 if

 \smallskip

  1. $A(v,w)\in{\Bbb Z}_p$ for all
 $v,w$ in $R$

  2. if $h\in \Q_p^2$ satisfies  $A(h,v)\in {\Bbb Z}_p$
   for all $v\in R$, then $h\in R$.

   \smallskip

   Vertices of the tree $\T(\Q_p^2)$ are self-dual lattices.
 Two vertices $R,S$ are connected by an edge if
 $$\text{volume of $R\cap S$}=\frac 1p \text{volume of $R$}$$

 It can be shown that  $\T(\Q_p^2)$   is the Bruhat--Tits tree
 $\T_p$. Obviously, the group ${\rm Sp}_2(\Q_p)$
 acts on our tree by automorphisms.

 \smallskip

 {\bf 2.8. Modular tree.}
 Consider the following standard picture from arbitrary
 textbook on complex analysis. Consider the  Lobachevsky
 plane $L:{\rm Im}\, z>0$ and the triangle $\Delta$ with three vertices
 $0$, $1$, $\infty$  on the absolute ${\rm Im}\, z=0$.
 Consider the reflections of $\Delta$  with respect to
 the sides of $\Delta$. We obtain 3 new triangles
 $\Delta_1$,  $\Delta_2$,  $\Delta_3$.
 Then we consider the reflections of $\Delta_j$ with respect
 to their sides etc. We obtain a tilling of
 $L$ by   infinite triangles
 (with vertices in rational points of the absolute
  ${\rm Im}\, z=0)$.

  Vertices of the modular tree are enumerated
   by the triangles of the tilling. Two vertices
   are connected by an edge, if the corresponding
   triangles have
   a common side.

   The group ${\rm SL}_2(\Z)$ acts on the modular tree
   in the obvious way.

 \smallskip

 {\bf 2.9 Tree of pants.} Let $R$ be a compact Riemann surface.
 Fix a collection $C_1$, \dots $C_k$  of closed mutually disjoint
 geodesics on $R$.
The universal covering of $R$
 is the Lobachevsky plane.

  The coverings of the cycles $C_j$ are geodesics on $L$.
  Thus we obtain the countable family of
  mutually disjoint geodesics on $L$. They divide $L$
  into the countable collection of domains.

  Now we construct a tree. Vertices of the tree are
  enumerated by the domains on $L$ obtained above.
   Two vertices are connected by
  an edge, if the corresponding domains have a common side.

  The fundamental group $\pi_1(R)$
  of the surface $R$ acts on this tree in the obvious way.

  \medskip

 {\large\bf 3. Hierarchomorphisms}

  \addtocounter{sec}{1}
 \setcounter{equation}{0}

 \medskip

 {\bf 3.1. Large group of hierarchomorphisms.}
 Consider a group $\Gamma$ acting on a simplicial (or metric) tree $J$.
 Consider
  a partition of $J$ into a finite collection of right subtrees
  $S_1$, \dots $S_k$,
 i.e., the subtrees $S_j$ are mutually disjoint, and
 $\Vert(J)=\bigcup\Vert(S_j)$.
Let
$$g_1:S_1\to J,\dots, g_k:S_k\to J$$
be a collection of embeddings such that

\smallskip

1) the subtrees $g_j(S_j)$ are mutually disjoint;

2) $\bigcup\Vert(g(S_j))=\Vert(J)$.

\smallskip

Thus we obtain the bijection
$$g=\{g_j,S_j\}:\Vert(J)\to\Vert(J)$$
given by
$$g(a)=g_j(a)\qquad\text{if}\quad a\in\Vert(S_j)$$
We call such maps  {\it hierarchomorphisms}, see Picture 2.
Denote the group of all such hierarchomorphisms
by $\Spher^\circ(J,\Gamma)$.

\smallskip

{\bf 3.2. Action of hierarchomorphisms on absolute.}
Consider a hierarchomorphism $g=\bigl\{g_j,S_j\bigr\}$.
Let $\omega\in \Abs(J)$. Let $a_1,a_2,\dots$ be a way
leading to $\omega$. For a sufficiently large $N$ and for  some
$S_j$, we have $a_N,a_{N+1},\dots\in S_j$. Hence
$g_j(a_N),g_j(a_{N+1}),\dots\in g_j(S_j)$ is a way leading
to some point
$$\nu\in\Abs\bigl((g_j(S_j)\bigr)\subset\Abs(J).$$

\2

We assume
$$\nu=g(\omega).$$

Fix a point $\xi\in\Vert(J)$.
Under the previous notation, consider the sequence
$$n_M=\rho(\xi,a_M)-\rho(\xi,g_j(a_M)).$$
This sequence becomes a constant after a sufficiently
large $M$.
We denote this constant (the {\it pseudoderivative})  by
$$n(g,\omega)=n_\xi(g,\omega).$$
The following statement is obvious.

\smallskip

{\sc Proposition 3.1.} {\it For $g,h\in\Spher^\circ(J,\Gamma)$,
 $\omega\in\Abs(J)$,}
\begin{equation}
n(g h,\omega)=n(h,\omega)+n(g,h\omega).
\end{equation}

{\bf 3.3. Small group of hierarchomorphisms.}
Denote by $\Spher(J,\Gamma)$ the group
of transformations of the absolute induced by
elements $g\in\Spher^\circ(J,\Gamma)$.

The kernel of the canonical map
  $$ \Spher^\circ(J,\Gamma) \to  \Spher(J,\Gamma)$$
  consists of finite permutations of the set
  $\Vert(J)$.

 Obviously,
 the pseudoderivative $n(g,\omega)$
 is well defined for   $g\in\Spher(J,\Gamma)$.

  \smallskip

 {\bf 3.4. A variant: planar hierarchomorphisms.}
 Assume a simplicial tree $J$ be planar
 (this means, that for each vertex $a$ we fix the
 cyclic order on the set of edges containing
 $a$; it is the case
 in some of our examples.
  Then also we have a canonical cyclic order
 on the absolute.

 Now we can consider the group of hierarchomorphisms
 that
 preserves the cyclic order on the absolute.

  \medskip

  {\large \bf 4. Hilbert spaces $\H(J)$}

   \addtocounter{sec}{1}
  \setcounter{equation}{0}

  \medskip

  {\bf 4.1. Definition.}
  Let $J$ be a metric tree, let $0<\lambda<1$.
  Denote by $\H(J)$ the real Hilbert space spanned
  by the formal vectors $e_a$, where $a$ ranges in $\Vert(J)$,
  with inner products given by
  \begin{equation}
  \langle e_a, e_b\rangle=\lambda^{\rho(a,b)},\qquad
  \forall a,b\in\Vert(J).
  \end{equation}

  We must show that a system of vectors
  with inner products (4.1) can be realized in a Hilbert space.

  \smallskip

  {\bf 4.2. Existence of $\H(J)$.}
  Let $a$ be a vertex of $J$. Let $b_1,b_2,\dots$ be  the
 vertices adjacent  to $a$. Consider an arbitrary unit vector  $e_a$
  in a real infinite dimensional Hilbert
   space $\cH$.
   Consider a collection
  $L_{b_1}$, $L_{b_2}$,\dots of pairwise perpendicular
  two-dimensional planes%
\footnote{Subspaces $M_1$, $M_2$ in a Hilbert space
  are perpendicular iff there is an orthogonal system of vectors
  $u_1,u_2,\dots$, $v_1,v_2,\dots$, $w_1,w_2,\dots$ such that
  $M_1$ is spanned by the vectors $u_i$, $v_j$, and $M_2$ is spanned by
 the vectors $w_n$, $v_j$.}
     containing $e_a$.
  For each plane $L_{b_k}$,
 we draw a vector $e_{b_k}\in L_{b_k}$ such that
 $$\langle e_{b_k}, e_a\rangle=\lambda^{\rho(a,b_k)},$$
 see Picture 3.

 \3

  By the perpendicularity,
  $$\langle e_{b_k},e_{b_l}\rangle=
    \langle e_{b_k},e_{a}\rangle \cdot
    \langle e_a,e_{b_l}\rangle=\lambda^{\rho(b_k,b_l)}.$$

  Then we apply the following inductive process.
  Assume that for a subtree $S$ the
  required embedding  $\Vert(S)\to \cH$
  is constructed, i.e., we have a subspace $\H(S)\subset \cH$.
   Let $b\in\Vert(J)$, and $c\not\in\Vert(J)$
  be adjacent to $b$. Consider the two-dimensional plane
  $L_c\subset \cH$
  that
 contains $e_b$ and is perpendicular to   $\H(S)$.
 Let us draw a unit vector $e_c\in L_c$ such that
 $$\langle e_c,e_b\rangle=\lambda^{\rho(b,c)}.$$
 Thus we obtained the required embedding
 $\Vert(S)\bigcup \{b\}\to \cH$.

  There is a sufficient place in the Hilbert space, and thus we
  obtain the embedding $\Vert(J)\to \cH$.

  \smallskip

 {\sc Remark.} This geometric picture is especially pleasant,
 if lengths of all edges are equal.

  \smallskip

  {\bf 4.3. More formal description of $\H(J)$.}
  Consider an {\it affine} real infinite dimensional Hilbert space $\cK$,
  i.e., a Hilbert space, where the origin of coordinates
  is not fixed. Denote by $\|\cdot\|$ the length in $\cK$.
  Consider a collection
   of points $N_a\in \cK$, where $a\in \Vert(J)$,
  such that

  \smallskip

   1) if $[a,b]$, $[c,d]$ are different edges of $J$, then
   $N_aN_b\bot N_cN_d$;

   \smallskip

   2) for $[a,b]\in\Edge(J)$,
   $$\|N_aN_b\|^2=\lambda\rho(a,b)\|.$$

   \smallskip

   The existence of such embedding is obvious.

   \4

   By the Pythagoras theorem,
   $$\|N_bN_c\|^2=\lambda\rho(b,c)\qquad \forall b,c\in \Vert(J).$$

   Now let us apply the following
    standard Fock--Schoenberg construction (\cite{Fok}, \cite{Sch}).
   For an affine Hilbert space $\cK$, there exists a
   linear Hilbert space
   ${\rm Exp}(\cK)$ and an embedding $\phi:\cK\to {\rm Exp}(\cK)$
   such that for all $X,Y\in \cK$
   $$\langle \phi(X,\phi(Y)\rangle=\exp(-\|XY\|^2).$$

   Fix any origin of the coordinates in $\cK$.
   We can assume that   ${\rm Exp}(\cK)$ is the direct sum of
   all symmetric powers of $\cK$
     $${\rm Exp}(\cK)=
     \R\oplus \cK\oplus S^2 \cK\oplus S^3\cK\oplus \dots,$$
     and
     $$
     \phi(X)= e^{-\|X\|^2}\Bigl[1\oplus\frac X{1!}
     \oplus \frac{X^{\otimes 2}}{2!}
      \oplus \frac{X^{\otimes 3}}{3!}\oplus \dots\Bigr].
      $$

    It remains to
     apply the Fock--Schoenberg construction
     to the space $\cK$ constructed above.
    The vectors $\phi(N_a)$ satisfy the relations (4.1).

    \smallskip

    {\sc Remark.} The spaces $\H$ associated with a tree are present
    in Olshansky's paper \cite{Ols2}. In a implicit form,
    they are present in \cite{Ism} (without a tree).

    \smallskip

    {\bf 4.4. Action of the group of hierarchomorphisms in $\H(J)$.}
    Let a group $\Gamma$ acts on $J$ by isometries.
    Then $\Gamma$ acts in $\H(J)$ by the orthogonal
    operators\footnote{An orthogonal
   operator is an  invertible operator
    in a real Hilbert space preserving
    the inner product} of the Hilbert space $\H(J)$ by
    the formula
    \begin{equation}
    U(g)e_a=e_{ga}.
    \end{equation}

    Let now $g\in\Spher^\circ(J,\Gamma)$ be a hierarchomorphism.
    Define  the operators $U(g)$ by the same formula (4.2).

     \smallskip

    {\sc Theorem 4.1.} a) {\it The operators $U(g)$ are well defined
    and bounded.}

     \smallskip

    b) {\it Each operator $U(g)$ can be represented in the form
    $U(g)=A(1+R)$, where $A$ is an orthogonal operator and
    $R$ is an operator of  finite rank.}

      \smallskip

      The theorem is proved below in 4.6.

      \smallskip

      {\bf 4.5. The subspaces $\H(S)$.}
      Let $S$ be a subtree in $J$. Denote by $\H(S)$
      the subspace in  $\H(J)$ generated by the
      vectors $e_c$, where $c\in\Vert(S)$.
      Denote by $P(S)$ the operator of projection
      $\H(J)\to\H(S)$.

  \smallskip

{\sc Lemma 4.2.}
 {\it Let $S_1$, $S_2$ be two disjoint subtrees in $J$.
Let $b\in\Vert(S_1)$, $c\in\Vert(S_2)$ be the nearest vertices
of the subtrees $S_1$, $S_2$.}

 \smallskip

  a) {\it The sum $\H(S_1)+\H(S_2)$ is a topological direct
  sum in $\H(J)$.}

 \smallskip

  b) {\it Let $Q:\H(S_1)\to\H(S_2)$
  be the restriction of the projection operator $P(S_2)$ to
  $\H(S_1)$. Then the image of $Q$ is the line spanned by $e_c$,
  and the kernel of $Q$ is the orthocomplement in $\H(S_1)$ to $e_b$.}

 \smallskip

 {\sc Proof.} Let $h_1\in \H(S_1)$,  $h_2\in \H(S_2)$
 be unit vectors. The both statements are corollaries
 of the following inequalities
 $$\langle h_1, h_2\rangle\le \langle h_1, e_c\rangle;\qquad
 \langle h_1, h_2\rangle \le   \langle e_b, h_2\rangle.$$

 {\bf 4.6. Proof of Theorem 4.1.}
 Let $g=\{g_j,S_j\}\in\Spher^\circ(J)$ be a hierarchomorphism.
 Without loss of generality (see Lemma 1.1), we can assume that
 $S_j$ are bushes or single-point sets.

  By Lemma 4.2, the decomposition
  $$\H(J)=\bigoplus_j \H(S_j)$$
  is a topological direct sum.

  Consider the bilinear form
  $$Q(h_1,h_2)=\langle U(g)h_1, U(g)h_2\rangle-  \langle h_1, h_2\rangle $$
  on $\H(J)\times \H(J)$.
  It is sufficient to prove that $Q$ is a bounded form
  on $\H(J)\times \H(J)$ and the rank of $Q$ is finite.

  The matrix of $Q$ in the basis $e_a$ is
  $$Q(e_a,e_b)= \langle e_{ga}, e_{gb}\rangle-
         \langle e_{a}, e_{b}\rangle =
         \lambda^{\rho(ga,gb)}-\lambda^{\rho(a,b)}$$
  The matrix $Q(e_a,e_b)$ has the natural block decomposition
  corresponding to the partition
  $$\Vert(J)=\bigcup \Vert(S_j)$$
  It is sufficient to prove that each
  block has  finite rank.

   Thus, let $a$ ranges in $S_i$, $b$ ranges in $S_j$.
   If $S_i$ is an one-point space, then the required statement is obvious.

   Thus, we assume
   that $S_i$, $S_j$ are bushes (see 1.1). Let $u_i$, $u_j$
   be their roots.

    If $S_i=S_j$, then $Q(e_a,e_b)$ is the identical zero.

    Thus, assume $S_i\ne S_j$. Then
    \begin{align*}
    \rho(a,b)&=\rho(a,u_i)+\rho(u_i,u_j)+\rho(u_j,b);\\
    \rho(ga,gb)&=\rho(ga,gu_i)+\rho(gu_i,gu_j)+\rho(gu_j,gb)= \\
               &=\rho(a,u_i)+\rho(gu_i,gu_j)+\rho(u_j,b)
     .\end{align*}
     Thus,
     \begin{align*}
     Q(e_a,e_b)=\bigl[\lambda^{\rho(gu_i,gu_j)}-
        \lambda^{\rho(u_i,u_j)}\bigr]
        \cdot \lambda^{\rho(a,u_i)}\cdot \lambda^{\rho(b,u_j)} =\\
       ={\rm const} \cdot \langle e_{u_i}, e_a\rangle
                  \cdot         \langle e_{u_j}, e_b\rangle
     .\end{align*}

     Thus we obtain that
       the bilinear form $Q$ on $\H(S_i)\times\H(S_j)$
     is given by the formula
     $$Q(h_1,h_2)={\rm const} \cdot
                 \langle e_{u_i}, h_1\rangle  \cdot
                  \langle e_{u_j}, h_2\rangle
                  $$
   Thus the form $Q$
    on $\H(S_i)\times\H(S_j)$ is of rank $\le 1$.
                  This finishes the proof.

\smallskip

{\bf 4.7. Remark. Spaces $\H$ associated with $\R$-trees.}
Let we have a countable family $J_1$, $J_2, \dots$
of metric trees and let we have isometric emeddings
$\iota_k:J_k\to J_{k+1}$:
$$
\dots
\stackrel{\iota_{k-1}}\longrightarrow J_k
\stackrel{\iota_{k}}\longrightarrow J_{k+1}
\stackrel{\iota_{k+1}}\longrightarrow J_{k+2}
\stackrel{\iota_{k+2}}\longrightarrow \dots
$$
Let $\bold J$ be the direct limit (the union) of $J_k$.
Such spaces are called {\it $\R$-trees}.\footnote{up to a minor
variation of terminology}

Obviously, we have the chain of inclusions
$$
\dots\subset\H(J_k)\subset \H(J_{k+1})
 \subset \H(J_{k+2})\subset \dots
$$

Denote the inductive limit of this chain by
$\H({\bold J})$. Thus the Hilbert space $\H$ survives
for $\R$-trees.
Nethertheless, the analogue of Theorem 4.1 is wrong.

    \medskip

    {\large\bf 5. Boundary spaces}

     \addtocounter{sec}{1}
    \setcounter{equation}{0}

    \medskip

 In this Section, we construct some spaces $\E$
 of 'distributions' on the absolute of
 a metric tree. These spaces can be considered
 as an analogue of the Sobolev spaces on the spheres.
For the Bruhat--Tits trees, the spaces $\E$ are well-known,
see \cite{Car}.
We also construct the operator $\H\to\E$ of restriction
of a "function on tree" to
the absolute.

\smallskip

{\it In this Section,
   $J$ is a locally finite perfect metric tree.}

   \smallskip

 {\bf 5.1. Balls in absolute.}
 Let $S$ be a branch of $J$.
A ball $B[S]\subset \Abs(J)$ is the absolute of
the branch $S$. If we
delete the root of the $S$ and all edges containing the root,
 then $S$ will be disintegrated into
the finite collection of branches
$S^{(1)}$, $S^{(2)}$\dots, $S^{(k)}$.
Hence the ball $B[S]$ admits the canonical partition
\begin{equation}
B[S]=B[S^{(1)}]\cup \dots\cup B[S^{(k)}]
.\end{equation}
into the  balls  $B[S^{(k)}]$.

We define the topology on $\Abs(J)$ by the assumption that
 all the balls $B[S]$ are open-and-closed subset in $\Abs[S]$.
Obviously, $\Abs(J)$ is a completely discontinuous compact set.

\smallskip

{\sc Remark.} Obviously, hierarchomorphisms locally
preserve hierarchy  of balls
 on the absolute\footnote{Firstly,
 this  hierarchy structure on $p$-adic manifolds
 was mentioned in Addendum
 in Serre's book \cite{Ser}.}.
 Obviously,
spheromorphisms  are homeomorphisms of the absolute.
But preserving of the  hierarchy of balls is
a very rigid  condition on a homeomorphism.

\smallskip

{\bf 5.2. New notation in the space $\H(J)$.}
Let us fix a vertex $\xi\in\Vert(J)$.
Let $a,b\in\Vert(J)$. Consider the way
$a_0=a, a_1,\dots, a_l=b$ connecting $a$, $b$.
Assume
$$\theta(a,b)=2\min\rho(\xi,a_j).$$

We emphasis that this function has sense also
if $a$ or $b$ are points of the absolute,
and the value $\theta(a,b)$
 is finite except the case $a=b\in\Abs(J)$.

For $a\in\Vert(J)$, consider the vector $f_a\in\H(J)$
given by
$$f_a=\lambda^{-\rho(\xi,a)}e_a.$$
Then
$$\langle f_a,f_b,\rangle=\lambda^{-\theta(a,b)}.$$

{\sc Remark.} Let $S$ be a subtree in $J$ containing $\xi$.
 For $c\in\Vert(J)$,
consider the nearest vertex $b\in \Vert(S)$.
Then the projection of $f_c$ to $\H(S)$ is $f_b$.

\smallskip

{\bf 5.3. Measures on $\Abs(J)$ and
compatible systems of measures on $\Vert(J)$.}
Let $R\subset J$ be a subtree. We say that $R$ is  {\it complete}
if any $a\in\Vert(R)$ satisfies
 one of two following conditions
(see Picture 5).

\5

\smallskip

 1. Any vertex $b$ of $J$ adjacent to $a$ is contained in $R$.

 2. Only one vertex of $J$ adjacent to $a$ is contained in $R$

\smallskip

 Let $\partial R$ denote the boundary of $R$, i.e., the set
 of all vertices of the second type.

 We also assume $\xi\in\Vert(R)\setminus \partial R$.

 \smallskip

 Consider a real-valued measure
 (charge) $\mu$ of {\it finite variation}
 on $\Abs(J)$. Recall that any measure $\mu$
 of finite variation admits the canonical
 representation
 $$\mu=\mu^+-\mu^-,$$
 where $\mu^\pm$ are nonnegative finite measures, and
 for some (noncanonical) Borel subset $U\subset\Abs$,
 $$\mu^-(U)=0;\qquad \mu^+(\Abs\setminus U)=0.$$
 The {\it variation} of the measure $\mu$ is
 $${\rm var}(\mu)=\mu^+(U)+\mu^-(\Abs\setminus U).$$

 For a complete subtree $R$, denote by $u_1$, $u_2$, \dots
 the points of $\partial R$. For any $u_k$, there exists a unique
 branch $S_{u_k}\subset J$ such that $u_k$
 is the root of $S_{u_k}$ and $\xi\notin S_{u_k}$.

 Consider the measure $\mu_R$ defined on the finite set $\partial R$
 by
 $$\mu_R(u_j)=\mu\bigl(B[S_{u_j}]\bigr).$$
 Consider also the vector
 $$\Psi[\mu|R]=\sum\limits_{u_j\in\partial R}
            \mu\bigl(B[S_{u_j}]\bigr)\,f_{u_j}.$$

 Let $R_2\supset R_1$ be complete subtrees. Then we have
 the obvious retraction
 $$\eta^{R_2}_{R_1}:\,\,\Vert(R_2)\to\Vert(R_1):$$
 if $a\in \Vert(R_2)$, then $\eta^{R_2}_{R_1}(a)$
 is the nearest vertex of $R_1$.

  \smallskip

 {\sc Lemma 5.1.}  a) {\it $\mu_{R_1}$ is the image of  $\mu_{R_2}$
under the retraction $\eta^{R_2}_{R_1}(a)$.}

   \smallskip

  b) {\it The vector $\Psi[\mu|R_1]$ is the projection of
  $\Psi[\mu|R_2]$ to the subspace $\H(R_1)$.
  In particular,}
  $$\|\Psi[\mu|R_1]\|\le\|   \Psi[\mu|R_2]\|.$$

   \smallskip

   {\sc Proof.} Assertion a) is obvious, and  assertion b)
   follows from the last remark from 5.2.  \hfill$\square$

   \smallskip

Conversely, consider a family of complete subtrees
$$R_1\subset R_2\subset R_3\subset\dots,$$
such that $\bigcup R_j=J$. Let for each $j$ we have a measure
$\nu_j$ on $\partial R_j$,
 and $\eta^{R_{j+1}}_{R_j}\nu_{j+1}=\nu_j$
 for all $j$.
If $\sup {\rm var}(\nu_j)<\infty$, then there exists
a unique measure $\nu$ on $\Abs$ such that $\nu_j=\nu_{R_j}$.

 \smallskip

{\bf 5.4. Boundary spaces  $\E\subset\H$.}
Let $R_1\subset R_2\subset \dots$ be a sequence of complete subtrees
in $J$, and $\bigcup R_k=J$ (the construction below
do not depend on  choice of the sequence).

Let $\mu$ be a measure of finite variation on $\Abs(J)$.
We say that $\mu$ belongs to the class $\E=\E(J)$
if
$$\lim_{j\to\infty}\|\Psi[\mu|R_j]\|_{\H} <\infty.$$

{\sc Proposition 5.2.}
 {\it For $\mu,\mu'\in \E$, the following statements
hold.}

 a) {\it There exists the following limit in the space $\H(J)$}
 \begin{align}
 \Psi[\mu]&:=\lim_{j\to\infty} \Psi[\mu|R_j].\\
 \text{b)}\qquad\,\,\, \qquad \qquad \qquad    \|\Psi[\mu]\|_{\H}&=
           \lim_{j\to\infty}\|\Psi[\mu|R_j]\|_{\H}.\\
\text{c)}\qquad\qquad \qquad  \langle \Psi[\mu], \Psi[\mu'] \rangle_{\H}&=
        \lim_{j\to\infty}
         \langle \Psi[\mu|R_j], \Psi[\mu'|R_j] \rangle_{\H}.
         \end{align}

         {\sc Proof.} All statements follow from Lemma 5.1.  \hfill$\square$

         \smallskip

 Thus we obtain the embedding $\E(J)\mapsto\H(J)$ given
 by $\Psi:\mu\to\Psi[\mu]$.
 We define the inner product
 in $\E(J)$  by
 $$\langle \mu_1,\mu_2\rangle_{\E(J)}:=
            \langle \Psi[\mu_1], \Psi[\mu_2] \rangle_{\H(J)}.$$

 Denote by $\EE\subset\H$ the image of the embedding $\Psi$.
Denote by $\ov\EE$ the closure of $\EE$ in $\H$,
 and also denote by $\ov \E$ the completion of the space
 $\E$ with respect to the norm (5.3).

 \smallskip

 {\bf 5.5. More direct description of $\E$.}
We can write formally
  \begin{align}
  \|\mu\|^2_{\E}=
      \iint_{\Abs\times\Abs}\lambda^{-\theta(\omega_1,\omega_2)}
      d\mu(\omega_1)\, d\mu(\omega_2) ;\\
\langle \mu_1,\mu_2\rangle_{\E}=
         \iint_{\Abs\times\Abs}\lambda^{-\theta(\omega_1,\omega_2)}
         d\mu_1(\omega_1)\, d\mu_2(\omega_2).
         \end{align}
 These integrals are very simple, since the integrand
  $\lambda^{-\theta(\omega_1,\omega_2)}$ has only countable
  set of values. Nevertheless, generally
  (even for the Bruhat--Tits tree $\T_2$)
  for $\mu_1,\mu_2\in\E$,  these
  integrals diverge as  Lebesgue integrals.

Our limit procedure is equivalent to the
Riemann improper integration in the following sense.
Consider a complete subtree $R\subset J$
such that $\xi\in R$.
Then $J\setminus R$ is a union of disjoint branches
$S_1,\dots, S_k$. Thus
$$\Abs(J)=B[S_1]\cup\dots\cup B[S_k].$$
Let us define the  Darboux sum
$${\cal S}_R(\mu_1,\mu_2)=
\sum\limits_{i,j}
\Bigl\{\min\limits_{\omega_1\in B[S_i],\,\,\,\omega_2\in B[S_j]}
   \lambda^{-\theta(\omega_1,\omega_2)}\Bigr\} \mu_1(B[S_i])\mu_2(B[S_j])
   .$$

{\sc Remark.}
If $i\ne j$, then the value
 $\lambda^{-\theta(\omega_1,\omega_2)}$
 is a constant on $B[S_i]\times B[S_j]$.

 \smallskip

We have seen, that
\begin{equation}
R_2\supset R_1\qquad \Rightarrow\qquad {\cal S}_{R_1}(\mu,\mu)
 \le {\cal S}_{R_2}(\mu,\mu)
 .\end{equation}
 Now we can define the integral (5.5)
as the limit of these Darboux
   sums under refinement of the partition.
   A measure $\mu$ is contained in $\E$ iff the Riemann
   integral (5.5) is finite.

   After this, we can define the inner product in $\E$
   as the Riemann improper integral (5.6).

Nevertheless, the space $\H$ was essentially
used in the justification of this construction,
since the convergence of Darboux sums and
positivity of the integral (5.5)
are not obvious.

    \smallskip

  {\bf 5.6. Non-emptiness of $\E$.}

  {\sc Theorem 5.3.}
  a) {\it There exists $\sigma$, which belongs to  $0\le\sigma\le1$,
  such that the space $\E$ is zero for $\lambda<\sigma$
  and the space $\E$ is not zero for $\lambda>\sigma$.}

    \smallskip

  b) {\it If lengths of edges of $J$ are bounded, then $\sigma<1$.}

    \smallskip

  c) {\it Let $J$ contain a subtree $I$ that is isomorphic to
  the Bruhat--Tits tree $\T_p$ as a simplicial tree, and
  lengths of all edges of $I$ are $\le \tau$. Then
  $\sigma\le 1/\sqrt[2\tau]{p}$.}

     \smallskip

  d) {\it Assume lengths of edges of $J$
  are bounded away from zero. Let the
  number $s(N)$ of $a\in\Vert(J)$, satisfying $d_\s(\xi,a)\le N$,
  has exponential growth, i.e., $s(N)\le \exp(\alpha N)$ for some constant
  $\alpha$. Then $\sigma>0$.}

     \smallskip

  The proof of the Theorem is contained below in 5.7--5.11

   \smallskip

  {\bf 5.7. Expansion of $\|\Psi\|^2$ into
   series with positive terms.}
  Let $R_0\subset R_1\subset R_2\subset\dots$
   be a sequence of complete subtrees
  in $J$, and $\bigcup R_m=J$. We say that the sequence
  $R_j$ is {\it incompressible}
  if

  $1^\circ$. $R_0$ consists of the vertex $\xi$;

  $2^\circ$. for each $m$, there exists $u\in\partial R_m$
  such that $\Vert(R_{m+1})\setminus\Vert( R_m)$
  consists of vertices adjacent to $u$.

  Fix a measure (charge) $\mu$ on $\Abs$.

  Obviously,
   $\Psi[\mu|R_0]=\mu(\Abs)f_\xi$,
   and hence
   $$  \|\Psi[\mu|R_0]\|^2=\mu(\Abs)^2.$$

  Let us evaluate
  $$z^{(m)}(\lambda)=\|\Psi[\mu|R_{m+1}]\|^2_{\H}-
                    \|\Psi[\mu|R_{m}]\|^2_{\H}.$$
  Let $u$ be the vertex defined in $2^\circ$.
  Let $v_1,\dots, v_n\in\partial R_{m+1}$
   be the vertices adjacent to $u$, see Picture 6.

   \6

  Let $\mu_{R_{m+1}}(v_k)=t_k$
  (these numbers can be negative), respectively
  $\mu_{R_m}(u)=t_1+ \dots + t_n$.
  It is readily seen that
  \begin{multline}
  z^{(m)}(\lambda)=
  \Bigl(\lambda^{-2\rho(\xi,u)}\sum\limits_{k\ne l} t_kt_l+
  \lambda^{-2\rho(\xi,u)} \sum\limits_k \lambda^{-2\rho(u,v_k)}t_k^2\Bigr)
  -\lambda^{-2\rho(\xi,u)} \bigl(\sum t_k\bigr)^2=\\
  =\lambda^{-2\rho(\xi,u)} \sum\limits_k (\lambda^{-2\rho(u,v_k)}-1)\,t_k^2
  .\end{multline}
  First, we observe that this expression is completely determined
  by the measure $\mu$ and the vertex $u$. The subtrees
  $R_m$, $R_{m+1}$ are nonessential. Hence it is natural
  to denote  $z^{(m)}(\lambda)$ by $z_u(\lambda)$.

   Thus,
   \begin{align}
     \|\Psi[\mu]\|^2&=
    \mu(\Abs)^2 + \sum\limits_{m=1}^\infty  z^{(m)}(\lambda)= \\=
   &
 \mu(\Abs)^2 +
      \sum\limits_{u\in\Vert(J),u\ne\xi}z_u(\lambda)
      .\end{align}

   We emphasis that

     a) all summands of this series are positive;

     b) all summands $z_u(\lambda)$ are decreasing
     functions on $\lambda$ for $0\le\lambda\le1$.

     \smallskip

  {\bf 5.8. Existence of $\sigma$.}
  The Statement a) of Theorem 5.4 follows from the last observation
  of previous subsection.

     \smallskip

  {\bf 5.9. Existence of $\E$.}
  It is sufficient to prove c), since b) is a corollary of c).
  Furthemore, it is sufficient to prove nontriviality of $\E(I)$
  for the subtree $I$. Denote by $R_k$ the subtree of $I$,
  consisting of all vertices $a\in I$ such that the simplicial
  distance $d_{\s}(\xi,a)\le k$.

    Consider the uniform measure $\mu_{R_k}$
    on $\partial R_k$, i.e., the measure of each point
    is $1/(p^{k-1}(p+1))$. Obviously, the measures $\mu_k$
    form a compatible system of the measures, denote by $\mu$
    the inverse limit of the measures $\mu_{R_k}$.

    Let us estimate
\begin{align*}
\Psi[\mu|R_k]\|^2&= \frac 1{(p+1)^2 p^{2(k-1)}}
\|\sum\limits_{a\in \partial R_k} f_a\|^2=\\
&= \frac 1{(p+1)^2 p^{2(k-1)}} \sum\limits_{a,b\in \partial R_k}
\lambda^{-2\rho(a,b)}\le
\\
&\le \frac 1{(p+1)^2 p^{2(k-1)}}
\sum\limits_{j=0}^k \lambda^{-2\tau j}\cdot\left\{
  \begin{matrix}
  \text{number of pairs  $(a,b)\in\partial R_k$}\\
  \text{ such that
  $d_{\s}(a,b)=2(k-j)$} \end{matrix}\right\} =\\
&\!\!\!\!\!\!\!\!\!\!\!\!\!\!\!\!\!\!\!\!\!\!\!\!\!
=\frac 1{(p+1)^2 p^{2(k-1)}}
 \Bigl[ (p+1)p^{2k-1} +
 \sum\limits_{j=1}^{k-1}\lambda^{-2\tau j} (p+1)p^{k-1} (p-1) p^{k-j-1}
 +\lambda^{-2k\tau}(p+1)p^{k-1}\Bigr]
 \le \\
&\le\sum\limits_{j=0}^k \lambda^{-2\tau j}p^{-j}
.\end{align*}
If $\lambda^{2\tau} p>1$, then these sums are uniformly bounded
in $k$; hence $\mu\in\E(I)\subset\E(J)$.

\smallskip

{\bf 5.10. Localization.}

\smallskip

{\sc Lemma 5.4.} {\it Let $\mu\in\E$, and let $B[S]\subset\Abs$ be a
ball. Let $\nu$ be the restriction of $\mu$ to
$B[S]$ {\rm(}i.e., $\nu(A)=\mu(A\bigcap B([S])$
for any Borel subset $A\subset\Abs${\rm)}.
Then $\nu\in\E$.}

\smallskip

{\sc Proof.} We can assume $\xi\notin S$. Denote by $v$
the root of the branch $S$.
 The quantity $\|\mu\|_{\E}^2$ is
 the sum of the series $\sum z_u(\lambda)$ given by (5.8), (5.10).
 The series for $\|\nu\|_{\E}^2$ is  obtained from
 the series fo  $\|\mu\|_{\E}^2$  by the following
 operations.

1) For $u$ lying between $\xi$ and  $v$,
the summands  $z_u(\lambda)$ are changed in
a non-predictible way.

2) For any $u\in S$, the summand  $z_u(\lambda)$
does not change.

3) All other summands become zero.

Obviously, the new series $z_u(\lambda)$ is convergent.
     \hfill $\square$

 \smallskip

 {\sc Remark.}
 Consider a Borel subset $U$ in the absolute.
 Let $\nu$ be the resriction of $\mu\in\E$
 to $U$. Generally, $\nu\notin\E$.
Also, generally, $\mu^\pm\notin\E$.

\smallskip

{\bf 5.11. Lower estimate of $\sigma$.}
By Lemma 5.4, if $\E\ne 0$,
then there exists  a measure $\mu\in\E$
such that $\mu(\Abs)\ne 0$. For definiteness, assume
$\mu(\Abs)=1$.

Let $\sigma$ be a lower bound for lengths of edges.
Consider a complete subtree $R\subset J$
defined by the condition $d_{\s}(\xi,a)\le N$.
Consider the measure $\mu_R$ on $\partial R$.
In notation 5.7,
$$\|\Psi[\mu]\|^2=1+\sum\limits_{v\in \Vert J,\,\, v\ne\xi} z_u\ge
  \sum\limits_{u\in\partial R} z_u\ge
$$
(by formula (5.8))
\begin{equation}
\ge \lambda^{-2N\sigma} (\lambda^{-2\sigma}-1)
\sum\limits_{u\in\partial R} \mu_R(u)^2
.\end{equation}
 The number of points of $\partial R$ is less than $\exp\{\alpha N\}$,
 where $\alpha$ is a constant. Furthemore,
  $\sum_{u\in\partial R} \mu_R(u)=1$,
  hence the last
 expression is larger than
  $$ \ge \lambda^{-2N\sigma} (\lambda^{-2\sigma}-1) \exp\{-\alpha N\}.$$
   For a sufficiently small $\lambda>0$, the last expression
   tends to $\infty$ as $N\to\infty$.
   and thus $\E=0$.

 \medskip

 {\large\bf 6. Action of group of hierarchomorphisms in $\E$}

   \addtocounter{sec}{1}
  \setcounter{equation}{0}

  \medskip
 Let $J$ be a perfect locally finite metric tree.
 Let $\H(J)\supset\ov\EE(J)\simeq \ov\E(J)$ be the same spaces as above.
 Let a  group $\Gamma$ act on $J$ by isometries.
 Let $\Spher^\circ(J,\Gamma)$, $\Spher(J,\Gamma)$
 be the corresponding hierarchomorphisms groups.
 The group  $\Spher^\circ(J,\Gamma)$ acts in $\H(J)$
 by the operators $U(g)$ given by (4.1).

 \smallskip

 {\bf 6.1. Action of hierarchomorphisms in $\E$.}

 {\sc Proposition 6.1} a) {\it The space $\ov\EE(J)\subset\H(J)$
 is invariant with respect to  $\Spher^\circ(J,\Gamma)$.}

  \smallskip

 b) {\it For $g\in  \Spher^\circ(J,\Gamma)$, the restriction
 of the operator $U(g)$ to $\ov\EE$
depends only on the corresponding element
 $\widetilde g\in \Spher(J,\Gamma)$.}

     \smallskip

 c) {\it The action of  $\Spher(J,\Gamma)$ in $\E(J)\simeq\EE(J)$
  is given by
 \begin{equation}
 T_\lambda(\widetilde g) \mu(\omega)=
 \lambda^{n(g,\omega)}\cdot\mu(g\omega),
  \qquad\text{where $g\in\Spher(J,\Gamma)$},
 \end{equation}
 where the pseudoderivative
  $n(g,\omega)=n(\widetilde g,\omega)$ of a hierarchomorphism
  on the absolute was defined in {\rm 3.2},
 and $\mu(g\omega)$ is the image of the measure $\mu$
 under the transformation $\omega\mapsto g\omega$.}

 \smallskip

 {\sc Remark.} For $g\in \Gamma$, the operator
 $T_\lambda(g)$ is unitary.

    \smallskip

 {\sc Proof.} Fix $g\in \Spher^\circ(J,\Gamma)$.
 Let $R_0\subset R_1\subset\dots$ be a incompressible sequence
 of complete subtrees as in 5.7, $\bigcup R_k=J$.
 Consider the sequence
 $g\cdot\partial R_1$, $g\cdot\partial R_2$, \dots .
 There exists $l$ such that for all  $k\ge l$
 $$g\cdot\partial R_k  =
    \partial T_k\qquad\text{where $T_k$ is a complete subtree}
    .$$
 Hence,
 $$U(g)\Psi[\mu|R_k]=\Psi[\nu|T_k].$$
 where $\nu$ is some measure on $\Abs(J)$.

 We must show that the numbers $\|\Psi[\nu|T_k]\|$
 are bounded. Consider the expansion of $\|\Psi[\mu]\|^2$
and $\|\Psi[\nu]\|^2$ into the series $\sum z^k(\lambda)$,
see (5.9), (5.10). The summands with  numbers $<l$ are essentially
different, but this do not influence on the convergence.
Other summands are rearranged and multiplied by
the factors $\lambda^{n(g,\omega)}$.

But  $\lambda^{n(g,\omega)}$ has only finite number of values and hence
the series $\sum z^k(\lambda)$ for the measure $\nu$ is also convergent.
Thus $\nu\in\E(J)$.

The statement a) is proved, the statement b) is obvious,
and the statement c) follows from the same considerations.

\smallskip

{\bf 6.2. Almost orthogonality.}

{\sc Theorem 6.2.} {\it Let $g\in\Spher(J,\Gamma)$.
The operators $T_\lambda(g)$
in $\E(J)$ given by
{\rm(6.1)} admit the representation  $T_\lambda(g)=A(1+Q)$,
where $A$ is an orthogonal operator and $Q$ is a
finite rank operator.}

This statement follows from Theorem 4.1.
This can  also be proved directly from the explicit formulas
(5.6), (6.1).

    Address (spring 2001):

     Erwin Schr\"odinger Institute for Mathematical
             Physics

    Boltzmanngasse, 9, Wien 1020, Austria

    \smallskip

    Permanent address:
     Institute of Theoretical and Experimental Physics,

     Bolshaya Cheremushkinskaya, 25

     Moscow 117259

     Russia

     \smallskip

    {\tt e-mail neretin@main.mccme.rssi.ru}

  \qquad\quad  {\tt neretin@gate.itep.ru}

  \end{document}